\documentclass{gtart}

\def\ifplaintex{\expandafter\ifx\csname documentclass\endcsname\relax}

\def\gtp{{\mathsurround=0pt\it $\cal G\mskip-2mu$eometry \&\ 
$\cal T\!\!$opology $\cal P\!$ublications}}  % GT publications

\def\recd{{\small Received:\qua\receiveddate\ifx\reviseddate\relax
\else\qquad Revised:\qua\reviseddate\fi\par}} 

\def\lognumber#1{\def\thelognumber{#1}}
\def\volumenumber#1{\def\thevolumenumber{#1}}
\def\volumeyear#1{\def\thevolumeyear{#1}}
\def\papernumber#1{\def\thepapernumber{#1}}
\def\pagenumbers#1#2{\def\startpage{#1}\def\finishpage{#2}}
\def\published#1{\def\publishdate{#1}}

\def\received#1{\def\receiveddate{#1}}

\def\accepted#1{\def\accepteddate{#1}}

\def\asciiaddress#1{\def\theasciiaddress{#1}}

\long\def\asciiabstract#1{\long\def\theasciiabstract{#1}}

\let\\\par\let\thelognumber\relax\let\thevolumenumber\relax
\let\thepapernumber\relax\let\thevolumeyear\relax\let\startpage\relax
\let\finishpage\relax\let\publishdate\relax\let\receiveddate\relax
\let\reviseddate\relax\let\accepteddate\relax\let\theasciititle\relax
\let\theasciiauthors\relax\let\theasciiaddress\relax
\let\theasciiabstract\relax

\let\theasciiemail\relax

\ifplaintex
\font\logobig=cmssbx10 scaled 3836
\font\logomed=cmssbx10 scaled 2557
\else
\font\logobig=cmssbx10 scaled 4200
\font\logomed=cmssbx10 scaled 2800
\fi

\long\def\makeagttitle{   %%% start of definition of \makeagttitle
\count0=\startpage
\agt\hfill      %   Journal title (top left) 
\hbox to 45truept{\vbox to 0pt{\vglue -13truept{\logomed A\kern -.37em{\logobig 
T}\kern -.38em G}\vss}\hss}
\break
{\small Volume \thevolumenumber\ (\thevolumeyear)
\startpage--\finishpage\nl
Published: \publishdate}

\vglue .25truein

{\parskip=0pt\leftskip 0pt plus
1fil\def\\{\par\smallskip}{\Large\bf\thetitle}\par\medskip} \vglue
0.05truein

{\parskip=0pt\leftskip 0pt plus 1fil\def\\{\par}{\sc\theauthors}
\par\medskip}%
 
\vglue 0.03truein 

{\small\leftskip 25truept\rightskip 25truept{\bf Abstract}\stdspace\theabstract

{\bf AMS Classification}\stdspace\theprimaryclass
\ifx\thesecondaryclass\relax\else; \thesecondaryclass\fi\par
{\bf Keywords}\stdspace \thekeywords\par}\vglue 7truept

}   %%%% end of definition of \makeagttitle

\ifplaintex
\font\phead=cmsl9 scaled 950
\font\pnum=cmbx10 scaled 913
\font\pfoot=cmsl9 scaled 950
\headline{\vbox to 0pt{\vskip -4.5mm\line{\small\phead\ifnum
\count0=\startpage ISSN 1472-2739 (on-line) 1472-2747 (printed)
\hfill {\pnum\folio}\else\ifodd\count0\def\\{ }% 
\ifx\theshorttitle\relax\thetitle\else\theshorttitle\fi\hfill{\pnum\folio}
\else\def\\{ and }{\pnum\folio}\hfill\ifx\theshortauthors\relax\theauthors
\else\theshortauthors\fi\fi\fi}\vss}}
\footline{\vbox to 0pt{\vglue 0mm\line{\small\pfoot\ifnum\count0=\startpage
\copyright\ \gtp\hfill\else
\agt, Volume \thevolumenumber\ (\thevolumeyear)\hfill\fi}\vss}}
\else
\font=cmsl9 scaled 1050
\font\lnum=cmbx10 
\font=cmsl9 scaled 1050
\makeatletter
\def\@oddhead{{\small\lhead\ifnum\count0=\startpage ISSN 1472-2739 
(on-line) 1472-2747 (printed)\hfill {\lnum\number\count0}\else\ifodd\count0
\def\\{ }\ifx\theshorttitle\relax \thetitle \else\theshorttitle\fi\hfill
{\lnum\number\count0}\else\def\\{ and }{\lnum\number\count0}
\hfill\ifx\theshortauthors\relax 
\theauthors\else\theshortauthors\fi\fi\fi}}\def\@evenhead{\@oddhead}
\def\@oddfoot{\small\lfoot\ifnum\count0=\startpage\copyright\ \gtp\hfill\else
\agt, Volume \thevolumenumber\ (\thevolumeyear)\hfill\fi}
\def\@evenfoot{\@oddfoot}
\makeatother
\fi
\let\maketitlepage\makeagttitle
\let\makeshorttitle\maketitlepage
\let\maketitle\maketitlepage

\newwrite\gtoutfile
\long\gdef\makeheadfile{  %%% start of definition of \makeheadfile
{\def\\{, }\def\s{ }
\immediate\openout\gtoutfile head.xxx
\immediate\write\gtoutfile{To: math@arxiv.org}
\immediate\write\gtoutfile{Subject: put OR rep NNNNN:ppppp}
\immediate\write\gtoutfile{--text follows this line--}
\immediate\write\gtoutfile{Proxy-for: \ifx\theasciiauthors\relax
\theauthors\else\theasciiauthors\fi\s<\ifx\theasciiemail\relax\theemail\else\theasciiemail\fi>}
\immediate\write\gtoutfile{\noexpand\\}
\immediate\write\gtoutfile{Authors: \ifx\theasciiauthors\relax
\theauthors\else\theasciiauthors\fi}
{\def\\{ }\immediate\write\gtoutfile{Title: \ifx\theasciititle\relax
\thetitle\else\theasciititle\fi}}
\immediate\write\gtoutfile{Subj-class: GT or SG, GR etc}
\immediate\write\gtoutfile{MSC-class: \theprimaryclass\ifx\thesecondaryclass\relax\else, \thesecondaryclass\fi}
\immediate\write\gtoutfile{Journal-ref: Algebr. Geom. Topol. \thevolumenumber\s
(\thevolumeyear) \startpage-\finishpage}
\immediate\write\gtoutfile{Comments: Published by Algebraic and
Geometric Topology at}
\immediate\write\gtoutfile{\s\s\s  http://www.maths.warwick.ac.uk/agt/AGTVol\thevolumenumber/agt-\thevolumenumber-\thepapernumber.abs.html}
\immediate\write\gtoutfile{\noexpand\\}
\immediate\write\gtoutfile{}
\ifx\theasciiabstract\relax
\immediate\write\gtoutfile{\theabstract}\else
\immediate\write\gtoutfile{\theasciiabstract}\fi
\immediate\write\gtoutfile{}
\immediate\write\gtoutfile{\noexpand\\}
\immediate\write\gtoutfile{}
\immediate\closeout\gtoutfile}}  %%% end of definition of \makeheadfile

\def\maketitlepage{\makeagttitle\makeheadfile}
\let\makeshorttitle\maketitlepage
\let\maketitle\maketitlepage

\def\ifplaintex{\expandafter\ifx\csname documentclass\endcsname\relax}

\def\gtp{{\mathsurround=0pt\it $\cal G\mskip-2mu$eometry \&\ 
$\cal T\!\!$opology $\cal P\!$ublications}}  % GT publications

\def\recd{{\small Received:\qua\receiveddate\ifx\reviseddate\relax
\else\qquad Revised:\qua\reviseddate\fi\par}} 

\def\lognumber#1{\def\thelognumber{#1}}
\def\volumenumber#1{\def\thevolumenumber{#1}}
\def\volumeyear#1{\def\thevolumeyear{#1}}
\def\papernumber#1{\def\thepapernumber{#1}}
\def\pagenumbers#1#2{\def\startpage{#1}\def\finishpage{#2}}
\def\published#1{\def\publishdate{#1}}

\def\received#1{\def\receiveddate{#1}}

\def\accepted#1{\def\accepteddate{#1}}

\def\asciiaddress#1{\def\theasciiaddress{#1}}

\long\def\asciiabstract#1{\long\def\theasciiabstract{#1}}

\let\\\par\let\thelognumber\relax\let\thevolumenumber\relax
\let\thepapernumber\relax\let\thevolumeyear\relax\let\startpage\relax
\let\finishpage\relax\let\publishdate\relax\let\receiveddate\relax
\let\reviseddate\relax\let\accepteddate\relax\let\theasciititle\relax
\let\theasciiauthors\relax\let\theasciiaddress\relax
\let\theasciiabstract\relax

\let\theasciiemail\relax

\ifplaintex
\font\logobig=cmssbx10 scaled 3836
\font\logomed=cmssbx10 scaled 2557
\else
\font\logobig=cmssbx10 scaled 4200
\font\logomed=cmssbx10 scaled 2800
\fi

\long\def\makeagttitle{   %%% start of definition of \makeagttitle
\count0=\startpage
\agt\hfill      %   Journal title (top left) 
\hbox to 45truept{\vbox to 0pt{\vglue -13truept{\logomed A\kern -.37em{\logobig 
T}\kern -.38em G}\vss}\hss}
\break
{\small Volume \thevolumenumber\ (\thevolumeyear)
\startpage--\finishpage\nl
Published: \publishdate}

\vglue .25truein

{\parskip=0pt\leftskip 0pt plus
1fil\def\\{\par\smallskip}{\Large\bf\thetitle}\par\medskip} \vglue
0.05truein

{\parskip=0pt\leftskip 0pt plus 1fil\def\\{\par}{\sc\theauthors}
\par\medskip}%
 
\vglue 0.03truein 

{\small\leftskip 25truept\rightskip 25truept{\bf Abstract}\stdspace\theabstract

{\bf AMS Classification}\stdspace\theprimaryclass
\ifx\thesecondaryclass\relax\else; \thesecondaryclass\fi\par
{\bf Keywords}\stdspace \thekeywords\par}\vglue 7truept

}   %%%% end of definition of \makeagttitle

\ifplaintex
\font\phead=cmsl9 scaled 950
\font\pnum=cmbx10 scaled 913
\font\pfoot=cmsl9 scaled 950
\headline{\vbox to 0pt{\vskip -4.5mm\line{\small\phead\ifnum
\count0=\startpage ISSN 1472-2739 (on-line) 1472-2747 (printed)
\hfill {\pnum\folio}\else\ifodd\count0\def\\{ }% 
\ifx\theshorttitle\relax\thetitle\else\theshorttitle\fi\hfill{\pnum\folio}
\else\def\\{ and }{\pnum\folio}\hfill\ifx\theshortauthors\relax\theauthors
\else\theshortauthors\fi\fi\fi}\vss}}
\footline{\vbox to 0pt{\vglue 0mm\line{\small\pfoot\ifnum\count0=\startpage
\copyright\ \gtp\hfill\else
\agt, Volume \thevolumenumber\ (\thevolumeyear)\hfill\fi}\vss}}
\else
\font=cmsl9 scaled 1050
\font\lnum=cmbx10 
\font=cmsl9 scaled 1050
\makeatletter
\def\@oddhead{{\small\lhead\ifnum\count0=\startpage ISSN 1472-2739 
(on-line) 1472-2747 (printed)\hfill {\lnum\number\count0}\else\ifodd\count0
\def\\{ }\ifx\theshorttitle\relax \thetitle \else\theshorttitle\fi\hfill
{\lnum\number\count0}\else\def\\{ and }{\lnum\number\count0}
\hfill\ifx\theshortauthors\relax 
\theauthors\else\theshortauthors\fi\fi\fi}}\def\@evenhead{\@oddhead}
\def\@oddfoot{\small\lfoot\ifnum\count0=\startpage\copyright\ \gtp\hfill\else
\agt, Volume \thevolumenumber\ (\thevolumeyear)\hfill\fi}
\def\@evenfoot{\@oddfoot}
\makeatother
\fi
\let\maketitlepage\makeagttitle
\let\makeshorttitle\maketitlepage
\let\maketitle\maketitlepage

\newwrite\gtoutfile
\long\gdef\makeheadfile{  %%% start of definition of \makeheadfile
{\def\\{, }\def\s{ }
\immediate\openout\gtoutfile head.xxx
\immediate\write\gtoutfile{To: math@arxiv.org}
\immediate\write\gtoutfile{Subject: put OR rep NNNNN:ppppp}
\immediate\write\gtoutfile{--text follows this line--}
\immediate\write\gtoutfile{Proxy-for: \ifx\theasciiauthors\relax
\theauthors\else\theasciiauthors\fi\s<\ifx\theasciiemail\relax\theemail\else\theasciiemail\fi>}
\immediate\write\gtoutfile{\noexpand\\}
\immediate\write\gtoutfile{Authors: \ifx\theasciiauthors\relax
\theauthors\else\theasciiauthors\fi}
{\def\\{ }\immediate\write\gtoutfile{Title: \ifx\theasciititle\relax
\thetitle\else\theasciititle\fi}}
\immediate\write\gtoutfile{Subj-class: GT or SG, GR etc}
\immediate\write\gtoutfile{MSC-class: \theprimaryclass\ifx\thesecondaryclass\relax\else, \thesecondaryclass\fi}
\immediate\write\gtoutfile{Journal-ref: Algebr. Geom. Topol. \thevolumenumber\s
(\thevolumeyear) \startpage-\finishpage}
\immediate\write\gtoutfile{Comments: Published by Algebraic and
Geometric Topology at}
\immediate\write\gtoutfile{\s\s\s  http://www.maths.warwick.ac.uk/agt/AGTVol\thevolumenumber/agt-\thevolumenumber-\thepapernumber.abs.html}
\immediate\write\gtoutfile{\noexpand\\}
\immediate\write\gtoutfile{}
\ifx\theasciiabstract\relax
\immediate\write\gtoutfile{\theabstract}\else
\immediate\write\gtoutfile{\theasciiabstract}\fi
\immediate\write\gtoutfile{}
\immediate\write\gtoutfile{\noexpand\\}
\immediate\write\gtoutfile{}
\immediate\closeout\gtoutfile}}  %%% end of definition of \makeheadfile

\def\maketitlepage{\makeagttitle\makeheadfile}
\let\makeshorttitle\maketitlepage
\let\maketitle\maketitlepage

\lognumber{7}
\volumenumber{1}
\volumeyear{2001}
\papernumber{7}
\published{2 March 2001}
\pagenumbers{143}{152}
\received{16 November 2000}
\accepted{28 February 2001}

\usepackage{amsmath, amssymb}
\usepackage[dvips]{graphicx}

\makeatletter
\def\thebibliography#1 {\@thebibliography@{Gab87}\small\parskip0pt % 
plus2pt\relax}
\makeatother

\title{Brunnian links are determined by their complements}
\authors{Brian Mangum\\Theodore Stanford}
\address{Barnard College, Columbia University\\Department
of Mathematics\\New York, NY 10027, USA}
\secondaddress{New Mexico State University\\Department
of Mathematical Sciences\\Las Cruces, NM 88003, USA}
\email{mangum@math.columbia.edu, stanford@nmsu.edu}
\asciiaddress{Barnard College, Columbia University\\Department
of Mathematics\\New York, NY 10027, USA\\New Mexico 
State University\\Department
of Mathematical Sciences\\Las Cruces, NM 88003, USA}

\newtheorem{theorem}{Theorem}

\newtheorem{corollary}[theorem]{Corollary}
\theoremstyle{definition}

\newtheorem{algorithm}[theorem]{Algorithm}

\begin{document}

\begin{abstract}

If \( L_1 \) and \( L_2 \) are two Brunnian links with all
pairwise linking numbers \( 0 \), then we show that \(L_1 \) and \( L_2 \)
are equivalent if and only if they have homeomorphic complements.  In
particular, this holds for all Brunnian links with at least three
components.  If \( L_1 \) is a Brunnian link with all pairwise linking
numbers \( 0 \), and the complement of \( L_2 \) is homeomorphic to the
complement of \( L_1 \), then we show that \( L_2 \) may be obtained from
\( L_1 \) by a sequence of twists around unknotted components.  Finally, we
show that for any positive integer $n$, an algorithm for detecting an
$n$--component unlink leads immediately to an algorithm for detecting
an unlink of any number of components.  This algorithmic
generalization is conceptually simple, but probably
computationally impractical.

\end{abstract}
\asciiabstract{If L_1 and L_2 are two Brunnian links with all pairwise 
linking numbers 0, then we show that L_1 and L_2 are equivalent if and
only if they have homeomorphic complements.  In particular, this holds
for all Brunnian links with at least three components.  If L_1 is a
Brunnian link with all pairwise linking numbers 0, and the complement
of L_2 is homeomorphic to the complement of L_1, then we show that L_2
may be obtained from L_1 by a sequence of twists around unknotted
components.  Finally, we show that for any positive integer n, an
algorithm for detecting an n-component unlink leads immediately to an
algorithm for detecting an unlink of any number of components.  This
algorithmic generalization is conceptually simple, but probably
computationally impractical.}


\maketitle

\addtocounter{section}{1}

All spaces and maps will be taken to be piecewise-linear or smooth.  A
link \( L \) is a closed one--dimensional submanifold of \( S^3 \) with
a finite number of connected components.  For convenience, we will
assume an arbitrary orientation on a link \( L \), and an
arbitrary ordering of its components.  Two links \( L_1 \) and \(
L_2 \) are equivalent if there exists a homeomorphism of pairs
from \( (S^3, L_1) \) to \( (S^3, L_2) \).  Such a homeomorphism
is required to preserve the orientation on \( S^3 \), but not
required to preserve the orientation or ordering on the \( L_i \).
If two links \( L_1 \) and \( L_2 \) are equivalent, then their
complements \( S^3 - L_1 \) and \( S^3 - L_2 \) are homeomorphic.
A basic question is when the converse of this statement holds. For
knots (links with a single component), this problem was solved by
Gordon and Luecke, who showed that two knots are equivalent if and
only if they have homeomorphic complements (\cite{GL89}).  For
links with more than one component, there are many examples where
this result fails.  A simple way to construct such examples is to
take any link \( L_1 \) with an unknotted component \( K \), cut
along the disk bounded by \( K \) (which usually intersects some
of the other components of \( L_1 \)), twist one or more times,
and reglue.  The resulting link \( L_2 \) will have a complement
homeomorphic to that of \( L_1 \), and in many easy cases, \( L_2
\) can be shown to be inequivalent to \( L_1 \).

Thus knots are determined by their complements, but links of more than one
component in general are not.  Three questions arise:

1) Can we characterize those links that are determined by their
complements?  That is, can we characterize those links $L_1$ such
that, if $S^3 - L_1$ is homeomorphic to $S^3 - L_2$, then $L_1$ is
equivalent to $L_2$?

2) For a given link (or a family of links) $L_1$, can we characterize the
links $L_2$ such that $S^3 - L_1$ is homeomorphic to $S^3 - L_2$?

3) There are many interesting families of links which arise naturally in
various contexts.  (Boundary links, homology boundary links, pure links,
Brunnian links; links of bounded genus, bridge number, braid index,
unknotting number, etc) Is there any such family $\cal L$ such that if
$L,L^\prime \in \cal L$ and $S^3 - L$ is homeomorphic to $S^3 - L^\prime$,
then $L$ is equivalent to $L^\prime$?

For knots, of course, the Gordon--Luecke result answers all these questions.
For links of more than one component, however, none of these questions have
been addressed in the literature as far as we know.

We will define a class of links called HTB links, which is almost the
same as the class of Brunnian links.  We will show that for HTB links,
the answer to Question 3 is yes (Theorem \ref{thm:main}, our main
result).  Regarding Question 2, for an HTB link $L$, we will give a
characterization of those $L^\prime$ with complement homeomorphic to
that of $L$ (Theorem \ref{thm:nonbrunnian}).  With regard to Question
1, it follows from our answer to Question 3 and the twisting
construction (mentioned above and discussed further below) that an HTB
link is never determined by its complement (among all links) unless it
is the trivial link or a knot.

Let \( K_1 \) and \( K_2 \) be two components of a link \( L \).  There is
a well-defined linking number which indicates the homology class of \( K_1
\) in \( S^3 - K_2 \).  A link \( L \) is homologically trivial if this
linking number is \( 0 \) for any pair of components of \( L \) (a property
that does not depend on the orientations of \( K_1 \) and \( K_2 \)).  A
link is trivial if the components bound pairwise disjoint disks.  An \( n
\)--component link is Brunnian if every \( (n-1) \)--component sublink is
trivial.  We consider an empty link to be trivial, and therefore we
consider a knot to be a 1--component Brunnian link.  A component of a link
is unknotted if it bounds a disk in \( S^3 \), which may intersect the
other components of the link.  We define an {\bf HTB link} to be a link
that is both homologically trivial and Brunnian.  The class of HTB links
includes all Brunnian links with more than two components.

In 1892, Brunn constructed links which become trivial when any single
component is deleted, but there did not exist at the time rigorous
methods to show that these links were in fact nontrivial.  It was
not until 1961 that it was shown by Debrunner \cite{Deb61} (who also
generalized Brunn's results) that nontrivial Brunnian links exist.
Since that time, Brunnian links have become a standard subclass
of links, and have been utilized and generalized in a number of
ways.

Let \( L \) be a link, and let \( X_L \) be \( S^3 \) with an open
tubular neighborhood $\eta(L)$ of \( L \) removed.  The boundary
components of \( X_L \) are two--dimensional tori, each torus
associated to some component of \( L \).  Each torus separates \(
S^3 \) into an inside (containing the associated component of \( L
\)) and an outside.  Each torus has a unique isotopy class of
essential simple closed curve called the meridian which is
homologically trivial inside of the torus.  Likewise, each torus
has a unique isotopy class of essential simple closed curve called
the longitude which is homologically trivial outside the torus.
The meridian and the longitude form a basis for $H_1$ of the
torus.

The isotopy class of any unoriented simple closed curve on the
torus is determined by a slope \( r = p/q \in \mathbb Q \cup \{1/0
\}\), where \( p \) indicates how many times the curve wraps
around in the meridian direction and \( q \) indicates how many
times it wraps around in the longitude direction.  (Note that the
specific orientations of the longitude and meridian are not
important, but for \( r \) to be well-defined they must together
be oriented compatibly with the orientation of \( S^3 \).) Using
the slope \( r \), we may perform Dehn surgery on a component \( K
\) of a link by cutting out a solid torus neighborhood of \( K \)
from $S^3$ and regluing it so that the curve of slope \( r \) now
bounds a disk in the solid torus. Notice that if the slope is \(
1/0 \) then performing Dehn surgery is simply removing a solid
torus and replacing it the same way.

If \( L \) is a link with \( n \) components, then we will often
write \( L \) with an \( n \)--tuple after it.  A slope \( p/q \)
in the \( i \)th position will indicate that \( p/q \) surgery has
been performed on the \( i \)th component, and a \( * \) will
indicate that no surgery has been performed.  Sometimes, this
notation will refer to the link formed by the components with a
$*$ in the manifold resulting from the surgeries, and sometimes it
will refer to the complement of this link.  The meaning will be
clear from the context.  So, for example, \( L = L(*,*, \dots *)
\), and \( L \) with the first component deleted is \( L (1/0,*,*,
\dots *) \), whereas $L(1/0, \dots 1/0) \cong S^3$.

If \( K \) is an unknotted component, then \( 1/q \) Dehn surgery
on $K$ gives the same result as cutting $S^3$ and the link along
the disk bounded by \( K \), twisting \( q \) times, regluing, and
deleting \( K \). In particular, the link resulting from \( 1/q \)
surgery on an unknotted component of \( L \) is a link in \( S^3
\) with one fewer component.

As in \cite{GL89}, we will approach the proof of our main result
by considering those Dehn surgeries on an HTB link \( L \) which
produce \( S^3 \).  We start with the following theorem which,
while no surprise to experts, does not seem to appear in the
literature.

\begin{theorem}
\label{thm:s3surgery}
If \( L \) is a link with \( n \) components and if \(
L(x_{1}, x_{2}, \dots x_{n}) \cong S^3 \) for every \( n \)--tuple
\( (x_1, x_2, \dots x_n) \in \{1/0, -1, 1\}^{n} \) then \( L \)
is a trivial link.
\end{theorem}

\begin{proof}
We proceed by induction on \( n \), the number of components of \( L \).
The theorem is true for knots in \( S^{3} \) since, by the
cyclic surgery theorem in \cite{CGLS}, $1$ and $-1$ surgery on a
nontrivial knot cannot both produce a manifold with trivial fundamental
group.  Now, assuming the theorem is true for all links of \( n-1 \) or
fewer components, we will show the result holds for any link with \( n \geq
2 \) components.  Let \( L = K_{1} \cup K_{2} \cup \dots \cup K_{n} \) be a
link in \( S^{3} \) such that \( L(x_{1},x_{2}, \dots x_{n}) \cong
S^3 \) for all $n$--tuples \( (x_1, x_2, \dots x_n) \in \{1/0, -1,
1\}^{n} \)

If \( X_{L} = S^{3} - \eta(L) \) is reducible, then there is some
\( S^{2} \subset X_{L} \) which does not bound a ball.  Since any
\( S^{2} \subset S^{3} \) separates and bounds a ball on each
side, the sphere in \( X_{L} \) must split \( L \) into two
nonempty sublinks \( L_{1} \) and \( L_{2} \).  The inductive
hypothesis implies that each of these sublinks is the trivial
link. Furthermore, they are separated by a sphere, so their union
\( L \) is also the trivial link.  In particular, $L$ is trivial
if some torus $T_{i}$ on $\partial X_{L}$ has a compressing disk
$D$ in $X_{L}$ because the frontier of a regular neighborhood of
$T_{i} \cup D$ would be a reducing sphere.  Therefore, we may
assume that $X_L$ is irreducible with incompressible boundary.

Each component of \( L \) must be unknotted because performing
$1/0$ surgery on any other component results in a link that, by
induction, is trivial.  Furthermore, each pair of distinct
components \( K_{i} \) and \( K_{j} \) have \( lk(K_{i},K_{j}) = 0
\).  Suppose that \( lk(K_{1},K_{2}) = m \). Let \( K = L(1, *,
1/0, 1/0, \dots 1/0) \), a knot in \( S^{3} \).  Then \(
K(1-m^{2}) = L(1, 1, 1/0, 1/0, \dots 1/0) \) which is homeomorphic
to $S^{3}$ by our assumption.  See, for example, \cite{Rol76} as a
reference for this modification of surgery instructions.  However,
\( H_{1}(K(1-m^{2})) \cong \mathbb Z/(1-m^{2})\mathbb Z \), so
$m=0$. Therefore, the link resulting from $1$ or $-1$ surgery on
any one component is a link in $S^3$, and the meridian--longitude
pairs of the remaining components are the same as the meridian--longitude
pairs for the corresponding sublink of $L$.

Consider the links $L(-1,*, \dots,*)$ and $L(1,*, \dots,*)$.  Each
of these are $(n-1)$--component links which satisfy the hypothesis
of the theorem, so they are trivial by induction.  Therefore, the
torus $T_2 = \partial \eta(K_2) \subset X_L$ becomes compressible
after both $1$ and $-1$ Dehn fillings on $T_1 \subset \partial
X_L$.  Since $\Delta(1,-1) = 2 > 1$ and since $X_L$ is irreducible with 
incompressible boundary, by Theorem 2.4.4 of \cite{CGLS}
(generalized by Wu in \cite{Wu92}) we see that $X_L$ is either
$T^2 \times I$ or a cable space, a solid torus with a regular
neighborhood of a \( (p,q), q \geq 2, \) cable knot of the core
curve of the torus removed.

In either case, $X_L$ is Seifert fibered. Burde and Murasugi
produced a complete list in \cite{BM70} of all links in \( S^{3}
\) with Seifert fibered complement.  Any such link with more than
one component has nontrivial linking number between some pair of
components.  This contradicts that every pair of components of $L$
are algebraically unlinked, so $X_L$ must in fact be reducible,
proving the theorem.
\end{proof}

\begin{theorem}
\label{thm:trivial}
Let \( L \) be an HTB link.  Suppose there exist slopes \( r_i = p_i/q_i
\), such that \( q_i \ne 0 \) for all \( i \), and such that \( L (r_1,
r_2, \dots r_n) = S^3 \).  Then \( L \) is trivial.
\end{theorem}

\begin{proof}
When \( L \) has one component, this is the well-known result of Gordon and
Luecke.  Suppose \( L \) has \( n>1 \) components.  Because \( L \) is
homologically trivial, its \( i \)th meridian has order \( |p_i| \) in \(
H_1 (L (r_1, r_2, \dots r_n)) = H_1 (S^3) \).  Thus \( |p_i| = 1 \) for all
\( i \).  The \( n \)th component of \( L \) is unknotted because \( L \)
is Brunnian, and so \( L (*,* \dots *, 1/q_n) \) is a link in \( S^3 \).
The linking number of the \( i \)th and \( j \)th components of \( L (*,*
\dots *,1/q_n) \) is \( l_{ij} + q_n l_{in} l_{jn} \), where \( l_{ij} \),
\( l_{in} \), and \( l_{jn} \) are the linking numbers in \( L \).  (This
holds for any link in which \( q_n \) twists are made about an unknotted \(
n \)th component.) Thus \( L (*,* \dots *,1/q_n) \) is homologically
trivial.

Deleting the first component of \( L (*, \dots *,1/q_n) \) yields
\( L (1/0, *, \dots, 1/q_n) \), and since this link may
also be obtained by deleting the first component of \( L \) and
then performing Dehn surgery on the last component of the
resulting link, it is trivial.  The same holds for the \( i \)th
component, \( i<n \), and so \( L (*, \dots *, 1/q_n) \) is
Brunnian.  Thus, by induction, \( L (*, \dots *, 1/q_n) \) is
trivial, and therefore \( L (1/x_1, \dots 1/x_{n-1}, 1/q_n)
= S^3 \) for all integers \( x_1, \dots x_{n-1} \).  Gordon
and Luecke's result (\cite{GL89}) then implies that \( L (1/x_1,
\dots 1/x_{n-1}, *) \) is a trivial knot for all integers
\( x_1, \dots x_{n-1} \), and therefore \( L (1/x_1,
\dots 1/x_n) = S^3 \) for all integers \( x_1,
\dots x_n \).  By Theorem \ref{thm:s3surgery}, \( L \) is
trivial.
\end{proof}

Note that in the proof of Theorem~\ref{thm:trivial} 
we use a weaker version of
Theorem~\ref{thm:s3surgery}, namely
that the $x_i$ in the hypothesis
range over all possible $1/n_{i}$ slopes for any integers $n_{i}$, 
rather than just $0$, $1$, and $-1$.
In discussing a preliminary version of this work with Martin
Scharlemann, we also discovered a proof of Theorem
\ref{thm:trivial} that does not rely on Theorem
\ref{thm:s3surgery} at all.  Rather, the proof suggested by Scharlemann
uses a theorem of Gabai concerning Dehn fillings of
manifolds with taut foliations \cite{Gab87} to complete the
induction.

If the \( i \)th component of a link \( L \) is unknotted, then \( L (*,*,
\dots 1/q_i, *, \dots *) \) is the link obtained by twisting \( q_i \)
times around that component and then discarding it.  If \( L \) is a
nontrivial HTB link, then \( L (*,*, \dots 1/q_i, *, \dots *) \) cannot be
a trivial link (unless \( q_i = 0 \) ) because then \( L (1,1, \dots 1/q_i,
1, \dots 1) \) would be \( S^3 \), contradicting Theorem \ref{thm:trivial}.
This generalizes one result of Mathieu \cite {Mat92}, who showed that for a
two--component HTB link, twisting around one component always ties a
nontrivial knot in the other component.  (More generally, he determines
exactly when twisting around an unknotted disk in the complement of an
unknot can produce a nontrivial knot.) Mathieu's result together with that
of Gordon and Luecke implies Theorem \ref{thm:trivial} for two--component
links.

We can now prove our main result.

\begin{theorem}
\label{thm:main}
Let \( L_1 \) and \( L_2 \) be HTB links.  Then \( L_1 \) is equivalent to
\( L_2 \) if and only if \( S^3 - L_1 \) is homeomorphic to \( S^3 - L_2
\).
\end{theorem}

\begin{proof}
Let \( h\co S^3 - L_1 \to S^3 - L_2 \) be a homeomorphism. (By a
homeomorphism of link complements, we always mean an orientation
preserving homeomorphism, as is required for link equivalence.)
Since $L_{1}$ and $L_{2}$ are homologically trivial, the
longitudes of the components of each link are null-homologous in
the link complements.  These are the only isotopy classes of
essential simple closed curves on each boundary component of $S^3
- \eta(L_1)$ and $S^3 - \eta(L_2)$ that are null-homologous.
Since $h$ induces isomorphisms on the homology groups, $h$ maps
longitudes to longitudes.  The map induced by \( h \) on the
homology of the torus boundaries of the complement of tubular
neighborhoods of the \( L_i \) must be invertible, and therefore
the \( i \)th meridian of \( L_1 \) must map to a slope \( 1/q_i
\) of the \( i \)th component of \( L_2 \).

If every meridian of $L_{1}$ is taken to a meridian of $L_{2}$, then it is
clear from the definition of Dehn surgery that we can extend $h$ to a
homeomorphism $\hat h\co (S^{3},L_{1}) \to (S^{3},L_{2})$.  That is, the two
links are equivalent.

Assume without loss of generality that the first meridian of $L_{1}$ is
taken to some slope $1/q \ne 1/0$ of the first component of $L_{2}$.  Then
we can extend $h$ to a homeomorphism from $L_{1}(1/0,*,*...*)$ to
$L_{2}(1/q,*,*...*)$.  The former is trivial because the link $L_{1}$ is
Brunnian.  Thus, the latter is also trivial since it is a link in $S^{3}$
and the trivial link is determined by its complement.  Therefore,
$L_{2}(1/q,1,1,...1) \cong S^{3}$ so $L_{2}$ must be the trivial link by
Theorem \ref{thm:trivial}.  As above, the homeomorphism $h$ implies that \(
L_{1} \) is also the trivial link, proving that \( L_1 \) and \( L_2 \) are
equivalent.
\end{proof}

\begin{figure}[ht]
\begin{center}
{\scalebox{0.37}{\includegraphics{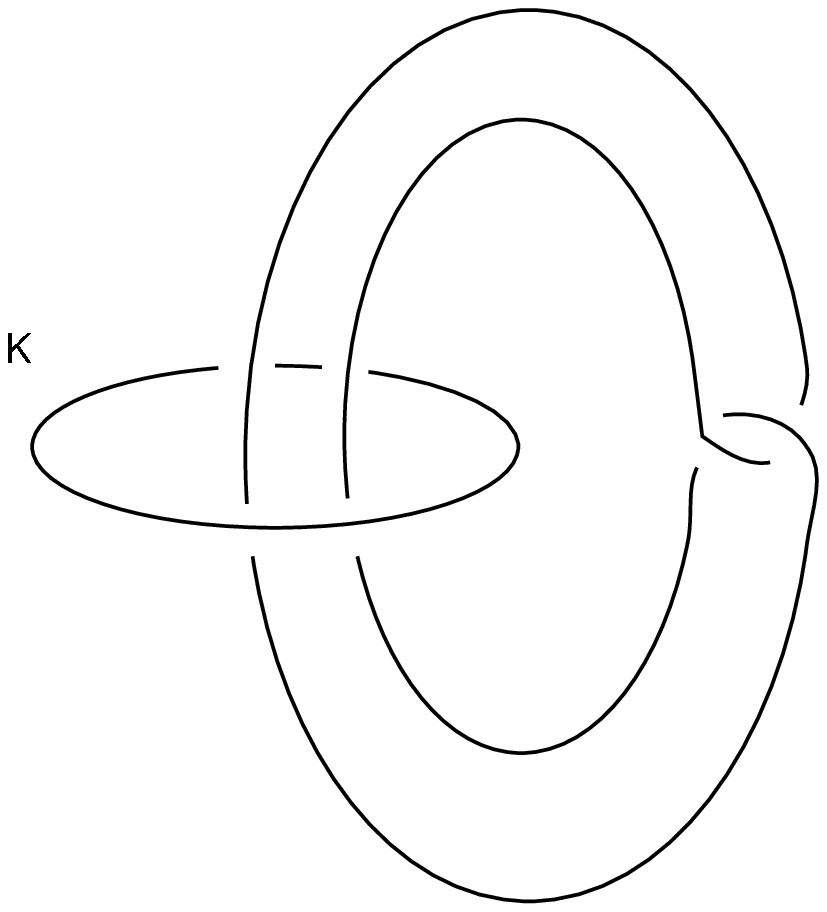}}}
$\ \ \longrightarrow \ \ $
{\scalebox{0.37}{\includegraphics{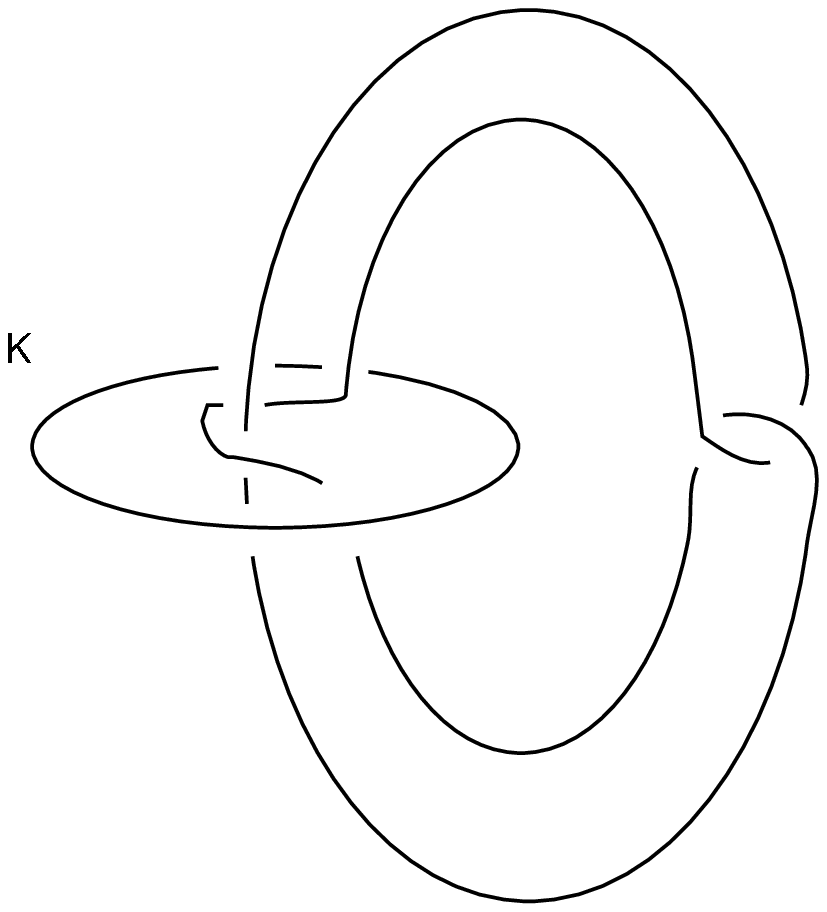}}}
\end{center}
\caption{\label{twisting} Inequivalent links formed by twisting along
unknotted $K$}
\end{figure}

As mentioned at the beginning of the paper, one way to construct
inequivalent links \( L_1 \) and \( L_2 \) with homeomorphic complements is
to twist one or more times around an unknotted component \( K \) of \( L_1
\) (without discarding \( K \)) to produce \( L_2 \).  (See Figure
\ref{twisting}.) More generally, a succession of such twisting operations
may be performed, as long as there remains an unknotted component to twist
around.  It follows from Theorem \ref{thm:main} that if \( L_1 \) is an HTB
link, then \( L_2 \) is never Brunnian and is therefore never equivalent to
\( L_1 \).  Moreover, if \( L_1 \) is an HTB link, then every \( L_2 \)
with a homeomorphic complement is obtained this way.  In contrast,
Berge (\cite{Ber91}) gives an explicit example of a pair of inequivalent 2
component links (not HTB links) with homeomorphic complements which remain
inequivalent even after twisting along the unknotted components.

\begin{theorem}
\label{thm:nonbrunnian}
Let $L_{1}$ be an HTB link and $L_{2}$ a link whose complement is
homeomorphic to that of $L_{1}$.  Then $L_{2}$ is obtained from $L_{1}$ by
successive twists around unknotted components.
\end{theorem}

\begin{proof}
Let $h\co S^3 - L_1 \to S^3 - L_2$ be a homeomorphism.  Since $L_{1}$ is
homologically trivial, the image of each longitude under $h$ must be an
essential simple closed curve on the boundary of $S^{3} - \eta(L_{2})$ that
is null-homologous in the link complement.  Therefore $L_{2}$ is also
homologically trivial, and $h$ maps longitudes to longitudes.

As in the proof of Theorem \ref{thm:main}, the \( i \)th meridian of \( L_1
\) must map to a slope \( 1/q_i \) of the \( i \)th component of \( L_2 \).
Moreover, the \( i \)th meridian of \( L_2 \) must then map under \( h^{-1}
\) to the slope \( -1/q_i \) of \( L_1 \).  This implies that \( L_2 \) is
determined by \( L_1 \) and the slopes \( q_i \), for if there is a
homeomorphism \( h^\prime\co S^3 - L_1 \to S^3 - L_2^\prime \), then the
composition \( h^\prime \circ h^{-1} \) takes the meridians of \( L_2 \) to
the meridians of \( L_2^\prime \), so \( L_2 \) and \( L_2^\prime \) are
equivalent.

Thus the links with complements homeomorphic to that of \( L_1 \) are
parameterized by \( n \)--tuples of slopes \( (1/q_1, 1/q_2, \dots 1/q_n)
\).  If \( q_i \ne 0 \) for all \( i \), then $L_1(-1/q_1, \dots, -1/q_n)
\cong S^3$ and by Theorem \ref{thm:trivial}, $L_1$ is trivial.  Therefore,
we assume that \( q_1 = 0 \).  Since \( L_1 \) is Brunnian, it has a
diagram where components \( 2,3, \dots n\) are \( n-1 \) disjoint planar
circles.  Using this diagram, it is clear that we may perform \( q_i \)
twists on the \( i \)th component for all \( 1 < i \le n \), and so in this
way, we may realize any link \( L_2 \) whose complement is homeomorphic to
\( S^3 - L_1 \).
\end{proof}

We now show how our results give ways of reducing the
problem of deciding if a given $n$--component link
is trivial to other decision problems.  First of all, recognizing the
trivial link may be reduced to recognizing \( S^3 \) itself, provided that
our method of representing links and three--manifolds allows Dehn surgery to
be done algorithmically.  Algorithms to recognize $S^{3}$ were produced by
Rubinstein in \cite{Rub95} and Thompson in \cite{Tho94}.

\begin{corollary}
\label{cor:trivial} A link \( L \) in \( S^3 \) with \( n \)
components is trivial if and only if \( L (e_1, e_2, \dots e_n) =
S^3 \) for each \( (e_1, e_2, \dots e_n) \in \{1/0,1\}^n \).
\end{corollary}

\begin{proof}
We proceed by induction on the number of components in $L$.  For
knots, the corollary is true by \cite{GL89}.  Suppose that $L$ has
more than one component.  It must be homologically trivial by the
argument in the proof of Theorem \ref{thm:s3surgery}.  Performing
$1/0$ surgery on any one component results in a trivial link by
induction, so $L$ is Brunnian.  Since $L(1, \dots, 1) \cong S^3$,
L is trivial by Theorem \ref{thm:trivial}.
\end{proof}

If we take a recursive approach, we can reduce the problem of recognizing a
trivial link to the problem of recognizing a trivial knot in \( S^3 \).
Suppose we have Algorithm A which tells whether or not a standard knot
diagram represents a trivial knot.  Then we may use the next corollary to
give an algorithm for recognizing a trivial link in \( S^3 \).

\begin{corollary}
\label{cor:trivial3}
A link \( L \) with \( n \) components is trivial if and only if
\begin{enumerate}
\item \( L \) is homologically trivial
\item Each \( (n-1) \)--component sublink \( L (*,* \dots 1/0, \dots *) \)
is trivial.
\item One of the \( (n-1) \)--component links \( L (*,*, \dots 1, \dots *)
\) is a trivial link in \( S^3 \).
\end{enumerate}
\end{corollary}

\begin{proof}
Conditions 1 and 2 imply that $L$ is an HTB link.  If condition 3
is true, then $L(1, \dots, 1) \cong S^3$ so L is trivial by
Theorem \ref{thm:trivial}.
\end{proof}

\begin{algorithm}
\label{alg}
An algorithm for recognizing a trivial link in \( S^{3} \).
\end{algorithm}

Given a diagram \( D \) for a link \( L \), the first step is to
calculate the pairwise linking numbers of \( L \).  This is simply
a matter of counting up positive and negative crossings in \( D
\).  If any linking number is nonzero, then \( L \) is not
trivial.  If all the pairwise linking numbers are zero, then pick
a component \( K \) of \( L \) and check whether it is trivial
using Algorithm A.  If \( K \) is not trivial, then \( L \) is not
trivial.  If \( K \) is trivial, then there exists a finite
sequence of Reidemeister moves which transform \( D \) into a
diagram \( D^\prime \) in which \( K \) is represented by a simple
closed curve with no self-crossings.  Given any link diagram,
there are only finitely many diagrams that result from a single
Reidemeister move on that diagram. Produce all diagrams that
result from a single Reidemeister move on $D$.  Then produce all
diagrams that result from a single Reidemeister move on any of
these, and so forth.  After finitely many steps, this will produce
the diagram \( D^\prime \).  
A twist around the disk bounded by \( K \) may then be
accomplished using the diagram \( D^\prime \). We now need to
check the triviality of \( n+1 \) \( (n-1) \)--component
links---namely the sublinks obtained by deleting \( 1 \) component
from \( L \) and the link \( L(1/1,*,*, \dots *) \). The algorithm
then proceeds recursively.

This algorithm is not likely to be of much practical use.  Even
if we had an efficient Algorithm A, we don't know a bound on the
number of Reidemeister moves it might take to convert a diagram
with an unknotted component $K$ into a diagram where $K$ is represented
by a circle.  Moreover, each twist along such a circle is likely
to raise the crossing number of the diagram substantially.

Suppose we have an $n$--unlink algorithm (ie, an algorithm for deciding
whether an $n$--component link is trivial).  Then we have an
$(n-1)$--unlink algorithm as well, by adding an extra disjoint and
unknotted component.  What we have seen is that by applying twisting
and deleting operations, we also have an $(n+1)$--unlink algorithm,
and thus our $n$--unlink algorithm becomes an unlink algorithm for
any number of components.

\bigskip

\noindent {\bf \large Acknowledgment:}

\smallskip

\noindent
The authors would like to thank Darren Long for very useful 
conversations regarding this work, and the referee for strengthening
Theorem~\ref{thm:s3surgery} and shortening its proof.
The second author was partially
supported by the US Naval Academy Research Council.

\Addresses\recd

\end{document}